# New existence results for prescribed mean curvature problem on balls under pinching conditions

Habib Fourti

*University of Monastir, Faculty of Sciences, 5019 Monastir, Tunisia.*
*Laboratory LR13ES21, University of Sfax, Faculty of Sciences, 3000 Sfax, Tunisia.*
*e-mail: Habib.Fourti@fsm.rnu.tn, habib40@hotmail.fr.*

**Abstract**

We consider a kind of Yamabe problem whose scalar curvature vanishes in the unit ball $\mathbb{B}^n$ and on the boundary $\mathbb{S}^{n-1}$ the mean curvature is prescribed. By combining critical points at infinity approach with Morse theory we obtain new existence results in higher dimensional case $n \geq 5$, under suitable pinching conditions on the mean curvature function.

## 1 Introduction and main results

This paper is concerned with the problem of prescribing mean curvature in the unit ball $\mathbb{B}^n$ of $\mathbb{R}^n$ with Euclidean metric $g_0$. Let $H$ be a smooth function on the boundary $\partial \mathbb{B}^n = \mathbb{S}^{n-1}$. Then, one may ask if there exists a conformal metric $g = u^{\frac{4}{n-2}} g_0$ whose scalar curvature vanishes in $\mathbb{B}^n$ and on the boundary $\partial \mathbb{B}^n$ the mean curvature is given by $H$. This geometric problem is equivalent to solving the following boundary value problem involving the Sobolev trace critical exponent

$$\begin{cases} \Delta u = 0 \quad \text{and} \quad u > 0 \quad \text{in } \mathbb{B}^n \\ \dfrac{\partial u}{\partial \nu} + \dfrac{n-2}{2} u = \dfrac{n-2}{2} H u^{\frac{n}{n-2}} \quad \text{on } \mathbb{S}^{n-1} \end{cases} \quad (1.1)$$

where $\nu$ is the outward unit normal vector on $\mathbb{S}^{n-1}$ with respect to the metric $g_0$.
The variational problem corresponding to (1.1) happens to be lacking of compactness, i.e. the functionals that we consider do not satisfy the Palais-Smale condition $(P.S)$ and this is due to the fact that the trace embedding $H^1(\Omega) \hookrightarrow L^{2(n-1)/(n-2)}(\partial \Omega)$ is not compact. This constitutes a strong obstruction for the application of the direct methods of the calculus of variations or even standard variational methods. Therefore, more refined techniques are needed as "Critical points at infinity theory" introduced by A. Bahri which we will follow in this work. The notion of critical point at infinity will be introduced by Definition 3.1 which is extracted from [7].
Besides the obvious necessary condition that $H$ must be positive somewhere, there is a Kazdan-Warner-type obstruction [22] to solve the problem, namely, if $u$ is a solution, then

$$\int_{\mathbb{S}^{n-1}} (\nabla_{g_0} H . \nabla_{g_0} x) u^{2(n-1)/(n-2)} d\sigma_{g_0} = 0,$$

where $x$ is the position vector of corresponding point of $\mathbb{S}^{n-1}$ in $\mathbb{R}^n$. Many works were devoted to the problem trying to understand under what conditions on $H$ equation (1.1) is solvable. Actually, several sufficient conditions were found (see the articles [1, 2, 3, 4, 12, 13, 15, 30, 31]). However, the question remained on whether they were or not necessary conditions. In this direction we present, here, a condition of new type to study the mean curvature problem on $\mathbb{S}^{n-1}$.

Before reviewing some known facts concerning problem (1.1) and stating our results, we point out that such a problem is a natural analogue of the problem of prescribing scalar curvature on the sphere known as the Nirenberg problem

$$(\mathcal{NP}) \quad -\Delta_g u + \frac{n(n-2)}{4} u = K u^{\frac{n+2}{n-2}}; \quad u > 0 \text{ in } \mathbb{S}^n,$$

where $\Delta_g$ denotes the Laplace Beltrami operator and $K$ is a smooth function on $\mathbb{S}^n$. Problem $(\mathcal{NP})$ has been investigated widely in the last decades, see for example [7, 8, 9, 11, 25, 27].

Due to Kazdan-Warner obstruction (see [24]), the solvability of the Nirenberg problem requires imposing conditions on the function $K$. It turns out that finding sufficient conditions under which the Nirenberg problem is solvable depends strongly on the dimension $n$ and the behavior of the function $K$ near its critical points. Indeed in low dimension $n < 5$ index counting criteria have been obtained, see [7, 9, 11, 25]. For dimension $n \geq 5$, such a counting index criterion fails, under the nondegeneracy assumption $(ND)$ (that is $\Delta K \neq 0$ at critical points of $K$). The extension to higher dimensions requires other conditions on the function $K$. We cite for instance [14] and [25] which deal with high dimensional spheres in respectively the perturbative setting (that is when $K$ is close to a constant) and the flat case. In particular, Chang-Yang obtained in [14] a perturbation theorem which asserts that there exists a conformal metric whose scalar curvature is equal to $K$, provided that the degree condition holds for $K$ which is a positive Morse function sufficiently close to $n(n-1)$ in $C^0$ norm.

A geometric flow has been introduced to study the Nirenberg problem by Struwe in [42] for $n = 2$, and has been generalized to $n \geq 3$ by Chen-Xu in [16]. Using the scalar curvature flow, Chen-Xu [16] was able to prove Chang-Yang's result with the quantitative bound on $\|K - n(n-1)\|_{C^0(\mathbb{S}^n)}$. Furthermore, regarding the high dimensional case $n \geq 5$, Malchiodi and Mayer [27] obtained recently an interesting existence criterion under some pinching condition. More precisely, they have extended the result in [16] to Einstein manifolds of positive Yamabe class under the pinching condition

$$\frac{K_{\max}}{K_{\min}} \leq 2^{\frac{1}{n-2}},$$

where $K_{\max} = \max_{\mathbb{S}^n} K$ and $K_{\min} = \min_{\mathbb{S}^n} K$. Moreover, under the following more stringent pinching requirement

$$\frac{K_{\max}}{K_{\min}} \leq \left(\frac{3}{2}\right)^{\frac{1}{n-2}},$$

and $K$ has at least two critical points with negative Laplacian they were able to prove that problem $(\mathcal{NP})$ has at least a solution. Very recently, Ahmedou and Ben Ayed [5] have given new existence results under pinching assumptions concerning a version of the Nirenberg problem on standard half spheres $\mathbb{S}^n_+$ with Neumann condition, which parallel the above mentioned existence results of Malchiodi-Mayer [27].

The aim of this paper is to extend these kind of results to the prescribing mean curvature equation. Before stating our results let us review some known facts. Going back to our problem (1.1), the special case $H = 1$ was considered by Escobar [19, 20, 21] and Beckner [10] who showed that the Yamabe ratio

$$S_n = \inf_{u \in H^1(\mathbb{B}^n)\setminus\{0\}} \frac{\int_{\mathbb{B}^n} |\nabla u|^2 dx + \frac{n-2}{2} \int_{\mathbb{S}^{n-1}} u^2 d\sigma_{g_0}}{\left(\int_{\mathbb{S}^{n-1}} u^{\frac{2(n-1)}{n-2}} d\sigma_{g_0}\right)^{\frac{n-2}{n-1}}}$$

is achieved by the conformal factors $u^{2/(n-2)} = |d\varphi|$ where $\varphi$ is a conformal transformation of the pair $(\mathbb{B}^n, \mathbb{S}^{n-1})$. More generally, Escobar[9], and later on, Li-Zhu [26] and Ou [29] have characterized such factors $u$ as the unique solutions of the equation (1.1) with $H \equiv 1$.

In [15], by a perturbation method, Chang, Xu and Yang obtained positive solutions by looking for constrained minimizers, more precisely, they proved that if $H$ is a smooth positive function on $\mathbb{S}^{n-1}$ and non-degenerate in the sense that $(\Delta_{\mathbb{S}^{n-1}} H)^2 + |\nabla_{\mathbb{S}^{n-1}} H|^2 \neq 0$ on $\mathbb{S}^{n-1}$, then under the degree condition:

$$\deg(G, \mathbb{B}^n, 0) \neq 1 \tag{1.2}$$

(1.1) has a solution for $\|H - 1\|_\infty$ sufficiently small. Here, $G : \mathbb{B}^n \to \mathbb{R}^n$ is some continuous map (see [15, p. 474] for the precise definition of $G$).

In [3, 4], the authors developed a Morse theoretical approach to this problem through a Morse-type reduction of the associated Euler-Lagrange functional in a neighborhood of its critical points at Infinity. This allowed them to provide some multiplicity and existence results under generic conditions on the function $H$. In fact, they presented in [4] the following index counting condition for the three dimensional case:

$$\sum_{y \in \mathcal{H}} (-1)^{ind(H,y)} \neq (-1)^{n-1}, \tag{1.3}$$

where

$$\mathcal{H} := \{y \in \mathbb{S}^{n-1} : \nabla H(y) = 0 \text{ and } \Delta H(y) < 0\} \tag{1.4}$$

and $ind(H, y)$ denotes the Morse index of $H$ at the critical point $y$.

Note that such a condition arises from computing the contribution of the degree of all critical points at infinity to the difference of topology between the level sets of the associated Euler-Lagrange functional. In dimension three, critical points at infinity with more than one mass were ruled out as we see in formula (1.3). Similarly, it was proven in [23] that blow-ups of solutions (of subcritical approximations) occur at most one point in the three dimensional case. However, in the four dimensional ball $\mathbb{B}^4$, the corresponding index counting condition (see [3] and [17]) is more complicated since critical points at infinity with more masses may occur.

Regarding the high dimensional case, Abdelhedi et al. have also provided in [4] a variety of classes of functions that can be realized as the mean curvature on the boundary of the standard $n$ dimensional ball, for $n \geq 5$, extending to problem (1.1) some topological invariants introduced by A. Bahri in his seminal work [8] concerning the scalar curvature problems in high dimensions.

To keep formula (1.3) in higher dimensions, the authors in [2] required that $H$ is flat near its critical points. Thanks to this flatness assumption they ruled out critical points at infinity with more than one bubble.

In addition to the previous conditions, by the method of geometric flow, Xu and Zhang presented in [31] a different Morse index one, known as an algebraic system condition, leading also to the solvability of problem (1.1). Moreover, they proved Chang-Xu-Yang's result [15] with the quantitative bound on $\|H - 1\|_{C^0(\mathbb{S}^{n-1})}$. To state their result let us introduce the following assumption.

(**A**): We assume that $H$ is a $C^3(\mathbb{S}^{n-1})$ positive function, which has only non-degenerate critical points with $\Delta H \neq 0$ (i.e. $H$ is enough smooth positive Morse function and non-degenerate on $\mathbb{S}^{n-1}$).

**Theorem 1.1** *[31] Let $n \geq 3$ and $H : \mathbb{S}^{n-1} \to \mathbb{R}$ satisfy the assumption (**A**) and the simple bubble condition*

$$\frac{\max_{\mathbb{S}^{n-1}} H}{\min_{\mathbb{S}^{n-1}} H} < \epsilon_n \tag{1.5}$$

*where $\epsilon_n = 2^{1/(n-1)}$, when $n = 3$ and $\epsilon_n = 2^{1/(n-2)}$, when $n \geq 4$. Let us consider the following numbers associated with $H$*

$$m_i = \sharp\{y \in \mathcal{H} : ind(H, y) = n - 1 - i\} \tag{1.6}$$

*where $ind(H, y)$ denotes the Morse index of $H$ at critical point $y$. If the following algebraic system has no non-trivial solutions,*

$$m_0 = 1 + k_0, \quad m_i = k_{i-1} + k_i, \quad 1 \leq i \leq n - 1, k_{n-1} = 0, \tag{1.7}$$

*with coefficients $k_i \geq 0$, then the problem (1.1) admits at least one positive solution.*

We recall some facts stated in [31]. Under the assumption that $H$ is a smooth positive Morse function on $\mathbb{S}^{n-1}$, the index counting condition (1.3) is equivalent to (1.2) following the appendix of Chang et al. [15]. On the other hand, one easily sees that (1.3) implies (1.7). If $H$ satisfies $\|H-1\|_{C^0(\mathbb{S}^{n-1})} < (\epsilon_n - 1)/(\epsilon_n + 1)$, then Theorem 1.1 implies that the boundary value equation (1.1) has a positive smooth solution and hence an estimate of $\|H-1\|_{C^0(\mathbb{S}^{n-1})}$ is obtained.

As we previously mentioned, the aim of the current paper is, by considering critical point at infinity theory combined with Morse theory, extending the result of Malchiodi and Mayer [27] to the case of mean curvature problem in high dimension $n \geq 5$. Precisely speaking, our first result can be stated as follows.

**Theorem 1.2** *Let $n \geq 5$ and $0 < H \in C^3(\mathbb{S}^{n-1})$ satisfying the assumption* (**A**). *If the following conditions hold*

1. $H_{\max}/H_{\min} < (3/2)^{1/(2(n-2))}$, *where* $H_{\max} := \max_{\mathbb{S}^{n-1}} H$ *and* $H_{\min} := \min_{\mathbb{S}^{n-1}} H$.

2. $\sharp \mathcal{H} \geq 2$, *where $\mathcal{H}$ is defined in (1.4) and $\sharp \mathcal{A}$ denotes the cardinality of the set $\mathcal{A}$.*

*Then the problem (1.1) has at least one solution.*

**Remark 1.3** *Compared to previous results that need one of the three topological conditions (1.2), (1.3) and (1.7) which are all sufficient for solvability of the boundary value problem (1.1), our result presents a condition of new type that is the assumption on the cardinality of critical points in $\mathcal{H}$. Moreover, our theorem is more practicable in the sense that, checking its conditions is easier than the ones of the previous results.*

The above pinching condition 1 of Theorem 1.2 can be relaxed when combined with some counting index formula. Namely, we prove

**Theorem 1.4** *Let $n \geq 5$ and $0 < H \in C^3(\mathbb{S}^{n-1})$ satisfying the assumption* (**A**). *If the following conditions hold*

1. $H_{\max}/H_{\min} < 2^{1/(2(n-2))}$.

2. $A_1 := \sum_{z \in \mathcal{H}} (-1)^{n-1-ind(H,z)} \neq 1$.

*Then the problem (1.1) has at least one solution.*

**Remark 1.5** *Although conditions (1.5) and (1.7) are weaker than respectively 1 and 2 of Theorem 1.4, this last theorem is very useful in our approach. Indeed, it is needed to prove our main result Theorem 1.2 which presents a condition of new type.*

**Remark 1.6**
- *We mention that pinching conditions in [5] and [27] are of order $((k+1)/k)^{1/(n-2)}$ where $k$ is some integer. However in our case they are of order $((k+1)/k)^{1/(2(n-2))}$. This little difference is due to the critical levels at infinity of our associated functional. See Proposition 3.4 and Corollary 3.6 which are the analogous results of [5, Corollary 3.18].*

- *Theorem 1.2 is similar to [5, Theorem 1.1] and our two main results are the analogous of [27, Theorem 1]. We point that the authors in [5] require that pinching condition is of order $(5/4)^{1/(n-2)}$. The argument developed in [5] requires to deal with critical point at infinity with less than or equal to four masses which justify the order of the pinching condition that is $(5/4)^{1/(n-2)}$ (obtained for $k = 4$). In our situation, we derive our main result from the analysis of the only two first critical levels at infinity and we obtain a pinching condition of second order ($k = 2$) reminiscent to the one given by [27].*

Next, we are going to describe our strategy for proving our main results. Our approach follows some arguments developed in [5] based on the techniques related to the critical points at infinity theory of Bahri [7] combined with Morse theory. As in [5], we start by describing the lack of compactness of the problem and characterizing the critical points at infinity of its associated functional which will be denoted by $J_H$. We collect these information from [4]. Then we compute the topological contribution of the critical points at infinity to the difference of topology between the level sets of the functional $J_H$.

Finally, two main ideas in [5] related to the pinching condition and inspired from [27] will lead to our existence results. The first one is that this kind of condition allows to show that suitable sublevels of $J_H$ are contractible. This will be the subject of our Proposition 4.2. Indeed, under such a condition, some energy levels of $J_H$ look very similar to those of the functional $J_{H\equiv 1}$ and since the only critical points of $J_{H\equiv 1}$ are minima lying in the bottom level $S_n$ and it has no critical points at infinity, we get all non empty sublevels of $J_{H\equiv 1}$ have trivial topology i.e. they are contractible sets. Hence the mentioned deformation lemma follows.

The second idea which is also related to the pinching condition is that critical levels at infinity of $J_H$ stratify depending on the number of bubbles.

Although the approach of Malchiodi and Mayer [27] is drastically different from ours involving a refined analysis for blowing up subcritical approximations, there is a big similarity between them and this is due to the narrow connection and correspondence between critical points at infinity of the associated functional and blow up solutions of the subcritical approximation problem.

The remainder of this paper is organized as follows. In section 2, we recall the variational framework and review the lack of compactness. The critical points at infinity, as well as their topological contributions are characterized in section 3 while section 4 is devoted to the proof of the main results.

## 2 The lack of compactness

Problem (1.1) has a variational structure, the Euler-Lagrange functional is

$$J_H(u) = \left( \int_{\mathbb{S}^{n-1}} H u^{\frac{2(n-1)}{n-2}} d\sigma_{g_0} \right)^{\frac{2-n}{n-1}}, \quad u \in H^1(\mathbb{B}^n).$$

The positive critical points of $J_H$, up to a multiplicative constant, are solutions of (1.1). Let $\Sigma$ be the unit sphere of $H^1(\mathbb{B}^n)$ equipped with the norm

$$\| u \|^2 = \int_{\mathbb{B}^n} |\nabla u|^2 dv_{g_0} + \frac{n-2}{2} \int_{\mathbb{S}^{n-1}} u^2 d\sigma_{g_0},$$

where $dv_{g_0}$ and $d\sigma_{g_0}$ denote the Riemannian measure on $\mathbb{B}^n$ and $\mathbb{S}^{n-1}$. We set $\Sigma^+ = \{u \in \Sigma / u > 0\}$. It is more convenient to work with the functional $J_H$ subjected to the constraint $u \in \Sigma$, since $J_H$ is lower bounded on $\Sigma$. Observe that the constrained critical points are indeed critical points of $J_H$. Thus, (1.1) is equivalent to finding the critical points of $J_H$ subjected to the constraint $u \in \Sigma^+$.

The exponent $\frac{2(n-1)}{n-2}$ is critical for the Sobolev embedding $H^1(\mathbb{B}^n) \hookrightarrow L^{\frac{2(n-1)}{n-2}}(\mathbb{S}^{n-1})$. This embedding being continuous and not compact, the functional $J_H$ does not satisfy the Palais-Smale condition, which leads to the failure of the standard critical point theory.

In order, to characterize the sequences failing the Palais-Smale condition, we need to introduce some notations. We will use the notation $x$ for the variables belonging to the unit ball $\mathbb{B}^n$ or to the half space $\mathbb{R}^n_+$ defined by $\mathbb{R}^n_+ := \{x \in \mathbb{R}^n, \ x_n > 0\}$. We will also use the notation $x = (x', x_n)$ for $x \in \mathbb{R}^n_+$. By a stereographic projection through an appropriate point in $\mathbb{S}^{n-1}$ we can reduce the problem to $\mathbb{R}^n_+$. Therefore, we will next identify the function $H$ and its composition with the stereographic projection $\pi$, and we will also identify a point $x \in \mathbb{B}^n$ by its image by $\pi$. See [4, p.333], for the expansion of $\pi$.

For $a \in \partial \mathbb{R}^n_+$ and $\lambda > 0$, we define the following function

$$\delta_{(a,\lambda)}(x) = \bar{c} \frac{\lambda^{\frac{n-2}{2}}}{((1+\lambda x_n)^2 + \lambda^2|x' - a'|^2)^{\frac{n-2}{2}}},$$

where $\bar{c}$ is chosen such that $\delta_{(a,\lambda)}$ is the family of solutions of the following problem

$$\begin{cases} \Delta u = 0 \quad \text{and} \quad u > 0 \quad \text{in } \mathbb{R}^n_+ \\ -\frac{\partial u}{\partial x_n} = u^{\frac{n}{n-2}} \quad \text{on } \partial \mathbb{R}^n_+. \end{cases}$$

Let $\tilde{\delta}_{(a,\lambda)}$ be the pull-back of $\delta_{(a,\lambda)}$ by the stereographic projection. We define now the set of potential critical points at infinity associated to the function $J_H$. For $\varepsilon > 0$ and $p \in \mathbb{N}$, let us define

$$V(p,\varepsilon) = \Big\{ u \in \Sigma \ / \ \exists a_1, ..., a_p \in \mathbb{S}^{n-1}, \exists \lambda_1, ..., \lambda_p > \varepsilon^{-1}, \exists \alpha_1, ..., \alpha_p > 0 \text{ with}$$
$$\| u - \sum_{i=1}^p \alpha_i \tilde{\delta}_{(a_i, \lambda_i)} \| < \varepsilon; \ | J_H(u)^{\frac{n-1}{n-2}} \alpha_i^{\frac{2}{n-2}} H(a_i) - 1 | < \varepsilon \ \forall i, \ \varepsilon_{ij} < \varepsilon \ \forall i \neq j \Big\},$$

where

$$\varepsilon_{ij} = \left( \frac{\lambda_i}{\lambda_j} + \frac{\lambda_j}{\lambda_i} + \frac{\lambda_i \lambda_j}{2} |a_i - a_j|^2 \right)^{-\frac{n-2}{2}}.$$

We are ready now to state the characterization of the sequences failing the Palais-Smale condition. Taking into account the uniqueness result of Li-Zhu [26] and using the idea introduced in [8] (see pages 325 and 334), we have the following proposition.

**Proposition 2.1** *Assume that (1.1) has no solution. Let $(u_k)$ be a sequence in $\Sigma^+$ such that $J_H(u_k)$ is bounded and $\partial J_H(u_k)$ goes to zero. Then there exist an integer $p \in \mathbb{N}$, a sequence $(\varepsilon_k) > 0$, $\varepsilon_k$ tends to zero, and an extracted subsequence of $u_k$'s, again denoted $(u_k)$, such that $u_k \in V(p, \varepsilon_k)$.*

We have the following proposition which defines a parametrization of the set $V(p, \varepsilon)$.

**Proposition 2.2** *For any $p \in \mathbb{N}$, there exists $\varepsilon_p > 0$ such that if $\varepsilon \leq \varepsilon_p$ and $u \in V(p, \varepsilon)$, then the following minimization problem*

$$\min_{\alpha_i > 0, \lambda_i > 0, a_i \in \mathbb{S}^{n-1}} \left\| u - \sum_{i=1}^p \alpha_i \tilde{\delta}_{(a_i,\lambda_i)} \right\|,$$

*has a unique solution $(\alpha, \lambda, a)$, up to a permutation.*

In particular, we can write $u$ as follows

$$u = \sum_{i=1}^p \alpha_i \tilde{\delta}_{(a_i, \lambda_i)} + v,$$

where $v$ belongs to $H^1(\mathbb{B}^n)$ and it satisfies the following condition:

$$< v, \psi > = 0 \quad \text{for } \psi \in \left\{ \tilde{\delta}_{(a_i,\lambda_i)}, \frac{\partial \tilde{\delta}_{(a_i,\lambda_i)}}{\partial \lambda_i}, \frac{\tilde{\delta}_{(a_i,\lambda_i)}}{\partial a_i} \right\} \quad i = 1, \ldots, p \quad (V_0)$$

here, $<.,.>$ denotes the scalar product defined on $H^1(\mathbb{B}^n)$ by

$$< u, v > = \int_{\mathbb{B}^n} \nabla u \nabla v \ dv_{g_0} + \frac{n-2}{2} \int_{\mathbb{S}^{n-1}} uv \ d\sigma_{g_0}, \quad u, v \in H^1(\mathbb{B}^n).$$

Here and in the next section we review some known facts extracted from [4] in order to characterize critical points at infinity associated to problem (1.1) under the assumption that (1.1) has no solution. We point out that, the authors in [4] have performed an analysis without such an assumption dealing with the possibility of existence of a critical point at infinity of new type consisting of sum of bubbles plus a solution $\omega$ of (1.1). All their expansions and results were done for $u - \alpha_0 \omega \in V(p, \varepsilon)$, where $\alpha_0 > 0$ satisfies $|\alpha_0 J_H(u)^{(n-1)/2} - 1| < \varepsilon$. They have also stated that their results hold true when $\omega = 0$.

We start by providing an expansion in $V(p, \varepsilon)$ of the functional $J_H$ on functions of the parameters $\alpha_i, \lambda_i, a_i, v$. This will allow us to understand the behavior of $J_H$ on the $v$-space.

**Proposition 2.3** [4] *For $\varepsilon > 0$ small enough and $u = \sum_{i=1}^{p} \alpha_i \tilde{\delta}_{(a_i, \lambda_i)} + v \in V(p, \varepsilon)$, we have the following expansion*

$$J_H(u) = \frac{S_n \sum_{i=1}^{p} \alpha_i^2}{\left(S_n \sum_{i=1}^{p} \alpha_i^{2\frac{n-1}{n-2}} H(a_i)\right)^{\frac{n-2}{n-1}}} \left[1 - \frac{n-2}{(n-1)\beta_1} \sum_{i=1}^{p} \alpha_i^{2\frac{n-1}{n-2}} \begin{cases} \frac{c_3 \Delta H(a_i)}{\lambda_i^2} \ln \lambda_i & \text{if } n = 3 \\ \frac{c_4 \Delta H(a_i)}{\lambda_i^2} & \text{if } n \geq 4 \end{cases}\right.$$

$$- \frac{c_1}{\gamma} \sum_{i \neq j} \alpha_i \alpha_j \varepsilon_{ij} + \frac{1}{\gamma}(Q_1(v, v) - f_1(v))$$

$$+ o\left(\sum_{i \neq j} \varepsilon_{ij} + \left(\sum_{i=1}^{p} \frac{\ln \lambda_i}{\lambda_i^2} \text{ if } n = 3; \sum_{i=1}^{p} \frac{1}{\lambda_i^2} \text{ if } n \geq 4\right)\right) + O\left(\|v\|^{\inf(3, 2\frac{n-1}{n-2})}\right) \right],$$

where,

$$f_1(v) = \frac{2\gamma}{\beta_1} \int_{\mathbb{S}^{n-1}} H\left(\sum_{i=1}^{p} \alpha_i \tilde{\delta}_i\right)^{\frac{n}{n-2}} v,$$

$$Q_1(v, v) = \|v\|^2 - \frac{n\gamma}{(n-2)\beta_1} \sum_{i=1}^{p} \int_{\mathbb{S}^{n-1}} H\left(\alpha_i \tilde{\delta}_i\right)^{\frac{2}{n-2}} v^2,$$

$$c_1 = \overline{c}^{\frac{2(n-1)}{n-2}} \int_{\mathbb{R}^{n-1}} \frac{dx}{(1+|x|^2)^{\frac{n}{2}}}, \quad c_3 = \overline{c}^4 \frac{\pi}{2}, \quad c_4 = \frac{\overline{c}^{\frac{2(n-1)}{n-2}}}{2(n-1)} \int_{\mathbb{R}^{n-1}} \frac{|x|^2}{(1+|x|^2)^{n-1}} dx$$

$$S_n = \overline{c}^{\frac{2(n-1)}{n-2}} \int_{\mathbb{R}^{n-1}} \frac{dx}{(1+|x|^2)^{n-1}}, \quad \gamma = S_n \sum_{i=1}^{p} \alpha_i^2 \quad \text{and} \quad \beta_1 = S_n \sum_{i=1}^{p} \alpha_i^{\frac{2(n-1)}{n-2}} H(a_i).$$

Then we deal with the $v$-part of $u$ to show that it can be neglected with respect to the concentration phenomenon.

**Proposition 2.4** [4] *For any $u = \sum_{i=1}^{p} \alpha_i \tilde{\delta}_{(a_i, \lambda_i)} \in V(p, \varepsilon)$, there exists a unique $\overline{v} = \overline{v}(\alpha, a, \lambda)$ which minimizes $J_H(u + v)$ with respect to $v \in H^1(\mathbb{B}^n)$ satisfying $(V_0)$. Moreover, we have the following estimates:*

$$\|\overline{v}\| \leq c \left[\sum_{i=1}^{p} \left(\frac{|\nabla H(a_i)|}{\lambda_i} + \frac{1}{\lambda_i^2}\right) + \sum_{i \neq j} \begin{cases} \varepsilon_{ij} \ln(\varepsilon_{ij}^{-1})^{\frac{1}{2}} & \text{if } n = 3, \\ \varepsilon_{ij}^{\frac{n}{2(n-2)}} \ln(\varepsilon_{ij}^{-1})^{\frac{n}{2(n-1)}} & \text{if } n \geq 4 \end{cases}\right].$$

Taking into account the Propositions 2.3 and 2.4, we end this section with an useful expansion of $J_H$ in $V(p, \varepsilon)$.

**Proposition 2.5** *[4] For the optimal $\overline{v}$ defined in the previous proposition, there exists a change of variables $v - \overline{v} \mapsto V$, such that,*

$$J_H(u) = \frac{S_n \sum_{i=1}^{p} \alpha_i^2}{(S_n \sum_{i=1}^{p} \alpha_i^{2\frac{n-1}{n-2}} H(a_i))^{\frac{n-2}{n-1}}} \left[ 1 - \frac{n-2}{(n-1)\beta_1} \sum_{i=1}^{p} \alpha_i^{2\frac{n-1}{n-2}} \begin{cases} \frac{c_3 \Delta H(a_i)}{\lambda_i^2} \ln \lambda_i & \text{if } n = 3 \\ \frac{c_4 \Delta H(a_i)}{\lambda_i^2} & \text{if } n \geq 4 \end{cases} \right.$$

$$\left. - \frac{c_1}{\gamma} \sum_{i \neq j} \alpha_i \alpha_j \varepsilon_{ij} + o\left( \sum_{i \neq j} \varepsilon_{ij} + \left( \sum_{i=1}^{p} \frac{\ln \lambda_i}{\lambda_i^2} \text{ if } n = 3;\ \sum_{i=1}^{p} \frac{1}{\lambda_i^2} \text{ if } n \geq 4 \right) \right) \right] + \| V \|^2.$$

## 3 Critical points at infinity and their topological contribution

Following A. Bahri, see [7, Definition 0.9], we set some useful definitions and notations.

**Definition 3.1** *A critical point at infinity of $J_H$ on $\Sigma^+$ is a limit of a non-compact flow line $u(s)$ of the equation:*

$$\begin{cases} \frac{\partial u}{\partial s} = -\partial J_H(u) \\ u(0) = u_0 \in H^1(\mathbb{B}^n). \end{cases}$$

*Assuming that (1.1) has no solution, $u(s)$ remains in $V(p, \varepsilon(s))$ for $s \geq s_0$.*
*Using Proposition 2.2, $u(s)$ can be written as:*

$$u(s) = \sum_{i=1}^{p} \alpha_i(s) \delta_{(a_i(s), \lambda_i(s))} + v(s).$$

*Denoting $a_i := \lim_{s \to +\infty} a_i(s)$ and $\alpha_i := \lim_{s \to +\infty} \alpha_i(s)$, we denote by*

$$(a_1, ..., a_p)_\infty \text{ or } \sum_{i=1}^{p} \alpha_i \delta_{(a_i, \infty)}$$

*such a critical point at infinity.*

As stated in Proposition 2.5, for $u = \sum_{i=1}^{p} \alpha_i \tilde{\delta}_{(a_i, \lambda_i)} + v \in V(p, \varepsilon)$ after change of variables we can write

$$J_H(u) = J_H\left( \sum_{i=1}^{p} \alpha_i \tilde{\delta}_{(a_i, \lambda_i)} + \overline{v} \right) + \|V\|^2.$$

We notice that in the $V$ variable, we define a pseudo-gradient by setting $\frac{\partial V}{\partial s} = -\mu V$, where $\mu$ is a very large constant. Then, at $s = 1$, $V(s) = e^{-\mu s} V(0)$, will be very small, as we wish. This shows that, in order to define our deformation, we can work as if $V$ was zero. The deformation will extend immediately, with the same properties, to a neighborhood of zero in the $V$ variable.

Using the estimates of the gradient vector field $\partial J_H$ (see Propositions 3.4 and 3.5 of [4]) and the expansion of $J_H$, we recall the characterization of the critical points at infinity associated with problem (1.1) which is based on the following key result.

**Proposition 3.2** *[4] Let $n \geq 5$. For any $u = \sum_{i=1}^{p} \alpha_i \tilde{\delta}_i \in V(p, \varepsilon)$ with $p \geq 1$ there exists a pseudo-gradient $W$ such that there exists a constant $c > 0$ independent of $u$ such that,*

1. $< -\partial J_H(u), W > \geq c \left( \sum_{i \neq j} \varepsilon_{ij} + \sum_{i=1}^{p} \left( \frac{|\nabla H(a_i)|}{\lambda_i} + \frac{1}{\lambda_i^2} \right) \right),$

2. $< -\partial J_H(u + \overline{v}), W + \frac{\partial \overline{v}}{\partial (\alpha_i, a_i, \lambda_i)}(W) >\geq c \left( \sum_{i \neq j} \varepsilon_{ij} + \sum_{i=1}^{p} \left( \frac{|\nabla H(a_i)|}{\lambda_i} + \frac{1}{\lambda_i^2} \right) \right).$

3. $W$ is bounded.

4. the only region where the maximum of the $\lambda_i$'s increases along the flow lines of $W$ is the region where $a_i$ is near a critical point $y_{j_i}$ of $H$ with, $-\Delta H(y_{j_i}) > 0$ and $j_i \neq j_r$ for $i \neq r$.

**Remark 3.3** *As stated in the end of the proof of [4, Proposition 4.1], we also have in the fourth statement that the $\lambda'_i s$ are of the same order.*

In the following, critical points at infinity in $V(p, \varepsilon)$ are characterized.

**Proposition 3.4** *[4] Let $n \geq 5$. Assume that $J_H$ does not have any critical point in $\Sigma^+$. Then, the only critical point at infinity of $J_H$ in $V(p, \varepsilon)$ for $\varepsilon$ small enough corresponds to $\sum_{j=1}^{p} \alpha_{i_j} \widetilde{\delta}_{(y_{i_j}, \infty)}$, where $\alpha_{i_j}^{-2} = S_n H(y_{i_j})^{n-2} \sum_{k=1}^{p} 1/H(y_{i_k})^{n-2}$, $y_{i_j}$ is a critical point of $H$ satisfying $-\Delta H(y_{i_j}) > 0$ and $i_j \neq i_k$ for $j \neq k$. Moreover, the Morse index of such critical points at infinity is equal to $p - 1 + \sum_{j=1}^{p} n - 1 - ind(H, y_{i_j})$, where $ind(H, y_{i_j})$ is the Morse index of $H$ at $y_{i_j}$.*

We introduce now a Morse Lemma at Infinity of $J_H$ near its critical points at infinity.

**Lemma 3.5** *Let $n \geq 5$ and $u = \sum_{i=1}^{p} \alpha_i \widetilde{\delta}_i + v \in V(p, \varepsilon)$ such that for each $1 \leq i \leq p$, $a_i$ is close to $y_i \in \mathcal{H}$ where $y_i \neq y_j$ for $i \neq j$. Then there exists a change of variables*

$$(\alpha_1, \ldots, \alpha_p, a_1, \ldots, a_p, \lambda_1, \ldots, \lambda_p, v) \mapsto \left(\widetilde{\alpha}, (A_1^+, A_1^-), \ldots, (A_p^+, A_p^-), \widetilde{\lambda}_1, \ldots, \widetilde{\lambda}_p, V\right),$$

*such that,*

$$J_H(u) = S_n^{\frac{1}{n-1}} \left(\sum_{i=1}^{p} \frac{1}{H(y_i)^{n-2}}\right)^{\frac{1}{n-1}} \left(1 - |\widetilde{\alpha}|^2 + \sum_{i=1}^{p} (|A_i^-|^2 - |A_i^+|^2)\right) \left(1 + \sum_{i=1}^{p} \frac{\widetilde{c}_i + \sigma}{\widetilde{\lambda}_i^2}\right) + \|V\|^2,$$

*where*

$$\widetilde{c}_i = -\frac{n-2}{n-1} \frac{c_4 \Delta H(y_i)}{S_n H(y_i)^{n-1} \sum_{j=1}^{p} \frac{1}{H(y_j)^{n-2}}},$$

*$\widetilde{\alpha} \in \mathbb{R}^{p-1}$ and $(A_i^+, A_i^-)$ are respectively the local coordinates of the parameters $\alpha(\overline{s}) := (\alpha_1(\overline{s}), \ldots, \alpha_p(\overline{s}))$ and $a_i(\overline{s})$. Here $\overline{s}$ is the real number defined in (3.8). Note that we lose an index for the parameter $\alpha$ since the functional $J_H$ is homogeneous with respect to this parameter.*

*Proof.* We argue as in [6, Lemma 5.4] whose proof is inspired from [8].
We recall that for $u = \sum_{i=1}^{p} \alpha_i \widetilde{\delta}_i + v$ there exists by Proposition 2.5 a change of variables such that

$$J_H(u) = J_H(\overline{u}) + \|V\|^2 \text{ where } \overline{u} := \sum_{i=1}^{p} \alpha_i \widetilde{\delta}_i + \overline{v}. \tag{3.1}$$

Since we have assumed that $y_i \neq y_j$, we get $|a_i - a_j| > 0$. Hence the parameter $\varepsilon_{ij}$ satisfies

$$\varepsilon_{ij} \leq c/(\lambda_i \lambda_i |a_i - a_j|^2)^{(n-2)/2} \leq c\left(\frac{1}{\lambda_i^{n-2}} + \frac{1}{\lambda_j^{n-2}}\right) = o\left(\frac{1}{\lambda_i^2} + \frac{1}{\lambda_j^2}\right) \text{ for } n \geq 5. \tag{3.2}$$

From Proposition 2.5 and by using the last estimate and the fact that $a_i$ is close to $y_i \in \mathcal{H}$ for each $i$, we obtain

$$J_H(u) = \frac{S_n \sum_{i=1}^{p} \alpha_i^2}{(S_n \sum_{i=1}^{p} \alpha_i^{2\frac{n-1}{n-2}} H(a_i))^{\frac{n-2}{n-1}}} \left[1 + \sum_{i=1}^{p} \frac{\widetilde{c}_i + o(1)}{\lambda_i^2}\right] + \|V\|^2.$$

Furthermore, for $\varepsilon' > 0$ small (with $\varepsilon/\varepsilon'$ small) and $W$ the pseudogradient constructed in Proposition 3.2, taken $\overline{u} = \sum_{i=1}^{p} \alpha_i \widetilde{\delta}_i + \overline{v} \in V(p, \varepsilon')$ we have:

$$< \partial J_H(\overline{u}), W > \leq -c \sum_{i=1}^{p} \left( \frac{|\nabla H(a_i)|}{\lambda_i} + \frac{1}{\lambda_i^2} \right), \tag{3.3}$$

where we have also used (3.2).

In the sequel, let $\sigma > 0$ be a small constant and set

$$I(\overline{u}) := \frac{S_n \sum_{i=1}^{p} \alpha_i^2}{(S_n \sum_{i=1}^{p} \alpha_i^{2\frac{n-1}{n-2}} H(a_i))^{\frac{n-2}{n-1}}} \left[ 1 + \sum_{i=1}^{p} \frac{\widetilde{c}_i + \sigma}{\lambda_i^2} \right]. \tag{3.4}$$

Since we have assumed that $1/\lambda_i$ is small for each $i$, it is easy to see that

$$0 < I(\overline{u}) - J_H(\overline{u}) < c\sigma \sum_{i=1}^{p} \frac{1}{\lambda_i^2}. \tag{3.5}$$

Moreover, it follows from Proposition 3.2 that

$$< \partial I(\overline{u}), W > \leq -c \sum_{i=1}^{p} \left( \frac{|\nabla H(a_i)|}{\lambda_i} + \frac{1}{\lambda_i^2} \right). \tag{3.6}$$

Next, for $\overline{u}_0 \in V(p, \varepsilon)$, we consider the following differential equation:

$$\frac{\partial u}{\partial s} = W(u); \quad u(0) = \overline{u}_0, \tag{3.7}$$

whose solution is $h_s(\overline{u}_0)$ where $h_s$ is the one-parameter group generated by $W$.

We note that, for $\overline{u} = \sum_{i=1}^{p} \alpha_i \widetilde{\delta}_i + \overline{v}$ as far as $h_s(\overline{u}) \in V(p, \varepsilon)$ such that $a_i$ is close to $y_i \in \mathcal{H}$ with $y_i \neq y_j$ for $i \neq j$, we have

$$h_s(\overline{u}) = \sum_{i=1}^{p} \alpha_i(s) \widetilde{\delta}_{a_i(s), \lambda_i(s)} + \overline{v(s)},$$

that is $\overline{v(s)}$ satisfies conclusions of Proposition 2.4.

We claim the following.

**Claim.** There exists $\overline{s} > 0$ such that

$$I(h_{\overline{s}}(\overline{u}_0)) = J_H(\overline{u}_0). \tag{3.8}$$

Indeed, since $I(h_s(\overline{u}_0))$ is a decreasing function with respect to $s$, there exists at most a unique solution to the equation $I(h_s(\overline{u}_0)) = J_H(\overline{u}_0)$. By using Remark 3.3, the only cases where there could be no solution are

- either $h_s(\overline{u}_0)$ exits $V(p, \varepsilon')$ (outside this set, we lose (3.3) since $W$ is only defined in $V(p, \varepsilon')$) before $I(h_{\overline{s}}(\overline{u}_0))$ reaches the level $J_H(\overline{u}_0)$.

- or $h_s(\overline{u}_0)$ will build a critical point at infinity before reaching this level.

In the following, we prove that these two cases cannot occur. In fact, for the first one, since $\overline{u}_0 \in V(p, \varepsilon)$ then, to exit $V(p, \varepsilon')$, the flow line has to travel from $V(p, \varepsilon'/2)$ to the boundary of $V(p, \varepsilon')$ by using again Remark 3.3. Note that, using (3.6), we have that: $\partial I(h_s(\overline{u}_0))/\partial s \leq -c(\varepsilon')$ along this path, independent of $\varepsilon$, but depending on $\varepsilon'$. Also, by (3.6), the time to travel from $V(p, \varepsilon'/2)$ to the boundary of $V(p, \varepsilon')$ is lower bounded by a constant $c'(\varepsilon')$ since $|W|$ is bounded and the distance to travel is lower bounded by a constant $c > 0$. Then $I(h_s(\overline{u}_0))$ decreases at least by $c(\varepsilon')c'(\varepsilon')$ during this trip.

However, using (3.5), since $\overline{u}_0 \in V(p, \varepsilon)$, it follows that $J_H(\overline{u}_0) < I(\overline{u}_0) \le J_H(\overline{u}_0) + c(\varepsilon)$. Hence, we have to choose $\varepsilon$ small with respect to $\varepsilon'$ so that $I(h_s(\overline{u}_0))$ crosses the level $J_H(\overline{u}_0)$ before $h_s(\overline{u}_0)$ leaves the set $V(p, \varepsilon')$.

Concerning the second case, the flow line $h_s(\overline{u}_0)$ will enter $V(p, \varepsilon_1)$ for each $\varepsilon_1 > 0$. Then, it follows from (3.5) that

$$0 < I(h_s(\overline{u}_0)) - J_H(h_s(\overline{u}_0)) \to 0 \text{ as } s \to \infty \text{ (for } \varepsilon_1 \to 0), \tag{3.9}$$

since $\lambda_i(s) \to \infty$ for each $i$ by using the fourth statement of Proposition 3.2 and Remark 3.3.
In the other hand, $J_H(h_s(\overline{u}_0))$ is a decreasing function of $s$ by using (3.3). Hence we get

$$J_H(h_s(\overline{u}_0)) < J_H(\overline{u}_0). \tag{3.10}$$

Through (3.9) and (3.10), $I(h_s(\overline{u}_0))$ has to cross the level $J_H(\overline{u}_0)$ and therefore this case cannot also happen.

Thereby, our claim is proved.

Conversely, taking $\varepsilon'' > 0$ small with respect to $\varepsilon$ and given $\overline{u}'_0 \in V(p, \varepsilon'')$, arguing by the same way (and using $-W$ as an increasing pseudogradient): there exists $\overline{s}' > 0$ such that $J_H(h_{-\overline{s}'}(\overline{u}'_0)) = I(\overline{u}'_0)$. Hence, $h_s$ is the required isomorphism.

It follows from (3.4) and the previous claim that

$$J_H(\sum_{i=1}^{p} \alpha_i \widetilde{\delta}_i + \overline{v}) = \frac{S_n \sum_{i=1}^{p} \alpha_i(\overline{s})^2}{\left(S_n \sum_{i=1}^{p} \alpha_i(\overline{s})^{2\frac{n-1}{n-2}} H(a_i(\overline{s}))\right)^{\frac{n-2}{n-1}}} \left[1 + \sum_{i=1}^{p} \frac{\widetilde{c}_i + \sigma}{\lambda_i(\overline{s})^2}\right].$$

Finally, it follows from classical Morse Lemma in finite dimensional case, that there exists a change of variables:

$$(\alpha_1(\overline{s}), \ldots, \alpha_p(\overline{s}), a_1(\overline{s}), \ldots, a_p(\overline{s})) \mapsto (\widetilde{\alpha}, (A_1^+, A_1^-), \ldots, (A_p^+, A_p^-))$$

such that $J_H$ reads as

$$J_H(\sum_{i=1}^{p} \alpha_i \widetilde{\delta}_i + \overline{v}) = S_n^{\frac{1}{n-1}} \left(\sum_{i=1}^{p} \frac{1}{H(y_i)^{n-2}}\right)^{\frac{1}{n-1}} \left(1 - |\widetilde{\alpha}|^2 + \sum_{i=1}^{p}(|A_i^-|^2 - |A_i^+|^2)\right) \left[1 + \sum_{i=1}^{p} \frac{\widetilde{c}_i + \sigma}{\widetilde{\lambda}_i^2}\right],$$

where $\widetilde{\lambda}_i = \lambda_i(\overline{s})$ for each $i$, $\widetilde{\alpha} \in \mathbb{R}^{p-1}$ and $(A_i^+, A_i^-)$ are respectively the local coordinates of the parameters $\alpha(\overline{s}) := (\alpha_1(\overline{s}), \ldots, \alpha_p(\overline{s}))$ and $a_i(\overline{s})$ (note that we lose an index for the parameter $\alpha$ since the functional $J_H$ is homogeneous with respect to this parameter). This ends the proof. □

Using the previous lemma, we identify the level sets of critical points at infinity.

**Corollary 3.6** *We will denote the corresponding critical point at infinity in Lemma 3.5 by $(y_1, \ldots, y_p)_\infty$. Such a critical point at infinity is at level*

$$C_\infty(y_1, \ldots, y_p) := S_n^{1/(n-1)} \left(\sum_{i=1}^{p} \frac{1}{H(y_i)^{n-2}}\right)^{1/(n-1)}.$$

We end this section with the topological contribution of the critical points at infinity to the difference of topology between the level sets of the functional $J_H$. As a consequence of the above corollary, Proposition 3.4 and the Morse reduction in Lemma 3.5 we have **Lemma 3.7** *Let $\tau_\infty$ be a critical point at infinity at the level $C_\infty(\tau_\infty)$ with index $i_\infty(\tau_\infty)$. Then for $\theta$ being a small positive number and a field $\mathbb{F}$, we have that*

$$H_l(J_H^{C_\infty(\tau_\infty)+\theta}, J_H^{C_\infty(\tau_\infty)-\theta}; \mathbb{F}) = \begin{cases} \mathbb{F} & \text{if} \quad l = i_\infty(\tau_\infty), \\ 0 & \text{otherwise,} \end{cases}$$

*where $H_l$ denotes the $l-$dimensional homology group with coefficients in the field $\mathbb{F}$.*

## 4 Proof of theorems

To prove Theorems 1.2 and 1.4 we follow the same strategy developed in [5] based on the characterization of the critical points at infinity in Proposition 3.4 and the computation of their contribution to the difference of topology in Lemma 3.7. But before using these two ingredients, we start by comparing some energy levels of $J_H$ to those of $J_{H\equiv 1}$, this will lead to the following deformation lemma. The idea behind is inspired from [27].

**Lemma 4.1** *Let $\underline{A}$ and $\overline{A} := (H_{\max}/H_{\min})^{(n-2)/(n-1)}\underline{A}$. Assume that $J_H$ does not have any critical point nor critical point at infinity in the set $J_H^{\overline{A}} \setminus J_H^{\underline{A}}$ where $J_H^A := \{u: J_H(u) < A\}$. Then for each $c \in (\underline{A}, \overline{A})$, the level set $J_H^c$ is contractible.*

*Proof.* Similar result is contained in [5] which deal with a Nirenberg type problem on $\mathbb{S}_+^n$ with Neumann condition. We follow the argument in [5, Lemma 4.1] which is inspired by the one developed in [27, Proposition 3.1].
Assuming that $J_H$ does not have any critical point nor critical point at infinity in $\Sigma^+$ between the levels $\overline{A}$ and $\underline{A}$, $J_H^{\overline{A}}$ retracts by deformation onto $J_H^{\underline{A}}$. Such a retraction can be realized by following the flow lines of a decreasing pseudogradient $W$ for $J_H$. In the sequel, we denote by $\phi_H$ the one parameter group corresponding to this pseudogradient and for each $u \in \Sigma^+$, we denote by $s_H(u)$ the first time such that $\phi_H(s_H(u), u) \in J_H^{\underline{A}}$.
Recall that, for $H \equiv 1$, the only critical points of $J_{H\equiv 1}$ are minima and lie in the bottom level $S_n$. Furthermore, in this case (when $H \equiv 1$), the Palais Smale condition holds along flow lines in $V(p,\varepsilon)$ for any $p \geq 1$ which assures that $J_{H\equiv 1}$ has no critical points at infinity. This result is similar to [7, Proposition 7.2 (i)] which deal with the Yamabe equation on the sphere.
In the sequel, we write $J_1$ instead of $J_{H\equiv 1}$. By following the flow lines of a decreasing pseudogradient $Z_1$ of the functional $J_1$, each flow line, starting from $u \in \Sigma^+$, will reach the bottom level $S_n$. Hence the set $J_1^A$ is a contractible one for each $A > S_n$. Let $\phi_1$ denote the one parameter group corresponding to $Z_1$.
Noticing that, we have

$$(1/H_{\max}^{(n-2)/(n-1)})J_1(u) \leq J_H(u) \leq (1/H_{\min}^{(n-2)/(n-1)})J_1(u) \text{ for each } u \in \Sigma^+,$$

we get

$$J_H^{\underline{A}} \subset J_1^{A'} \subset J_H^{\overline{A}} \text{ where } A' := H_{\max}^{(n-2)/(n-1)}\underline{A}.$$

Moreover we observe that for each $u \in \Sigma^+$, there exists a unique $s_1(u)$ satisfying $\phi_1(s_1(u), u) \in J_1^{A'}$. Next, we define the following map:

$$F := [0,1] \times J_1^{A'} \to J_1^{A'}; \quad F(t,u) := \phi_1(s_1(\phi_H(ts_H(u), u)), \phi_H(ts_H(u), u)).$$

We point out that $F$ is well defined and continuous and satisfies the following properties:

- For $t = 0$, we have $\phi_H(0, u) = u$. Furthermore, for each $u \in J_1^{A'}$, we have $s_1(u) = 0$. Therefore, for each $u \in J_1^{A'}$, we get $F(0, u) = \phi_1(0, u) = u$.

- For $t = 1$, we have $\phi_H(s_H(u), u) \in J_H^{\underline{A}} \subset J_1^{A'}$ (by the definition of $s_H$) which implies that $s_1(\phi_H(s_H(u), u)) = 0$ and therefore $F(1, u) = \phi_1(0, \phi_H(s_H(u), u)) = \phi_H(s_H(u), u) \in J_H^{\underline{A}}$ for each $u \in J_1^{A'}$.

- If $u \in J_H^A$, then $s_H(u) = 0$ which implies that $\phi_H(ts_H(u), u) = \phi_H(0, u) = u$. Therefore $F(t, u) = \phi_1(s_1(u), u) = \phi_1(0, u) = u$ for each $u \in J_H^A$ and each $t \in [0, 1]$ (we used $s_1(u) = 0$ since $u \in J_H^A \subset J_1^{A'}$).

Thus $J_1^{A'}$ retracts by deformation onto $J_H^A$, a fact which provides the claim of the lemma since $J_1^{A'}$ itself is a contractible set. □

We derive another deformation lemma from the previous one, the assumption (**A**) of this paper and an appropriate pinching condition imposed to the function $H$. Thanks to this last condition critical levels at infinity stratify (depending on the number of masses) leading to the contractibility of some energy sublevels of $J_H$. To introduce our second deformation lemma, we set for $\ell \in \mathbb{N}$

$$C_{\max}^{\ell,\infty} := (\ell S_n)^{1/(n-1)}/H_{\min}^{(n-2)/(n-1)} \quad \text{and} \quad C_{\min}^{\ell,\infty} := (\ell S_n)^{1/(n-1)}/H_{\max}^{(n-2)/(n-1)}.$$

Using Corollary 3.6, it is easy to see that the level of critical points at infinity corresponding to $\ell$ boundary points lies between $C_{\min}^{\ell,\infty}$ and $C_{\max}^{\ell,\infty}$. As [5, Proposition 4.2], we have the following:

**Proposition 4.2** *For $k \in \mathbb{N}$ being fixed, let $0 < H \in C^3(\mathbb{S}^{n-1})$ satisfying the condition (**A**) and the pinching condition $H_{\max}/H_{\min} < ((k+1)/k)^{1/(2(n-2))}$.*
*If $J_H$ does not have any critical point under the level $C_{\min}^{k+1,\infty}$. Then, for every $1 \leq \ell \leq k$ and every $c \in (C_{\max}^{\ell,\infty}, C_{\min}^{\ell+1,\infty})$, the sublevel $J_H^c$ is a contractible set.*

*Proof.* Assuming that $H_{\max}/H_{\min} < ((k+1)/k)^{1/(2(n-2))}$, we get for each $1 \leq \ell \leq k$, $(k+1)/k \leq (\ell+1)/\ell$ and

$$C_{\max}^{\ell,\infty} < C_{\max}^{\ell,\infty}(H_{\max}/H_{\min})^{(n-2)/(n-1)} < C_{\min}^{\ell+1,\infty}.$$

Indeed, we have

$$\begin{aligned}
C_{\max}^{\ell,\infty} &< C_{\max}^{\ell,\infty}(H_{\max}/H_{\min})^{(n-2)/(n-1)} \\
&< (\ell S_n)^{1/(n-1)}/H_{\min}^{(n-2)/(n-1)}((\ell+1)/\ell)^{\frac{1}{2(n-1)}} \\
&< ((\ell+1)S_n)^{1/(n-1)}/H_{\min}^{(n-2)/(n-1)}(\ell/(\ell+1))^{\frac{1}{2(n-1)}} \\
&< ((\ell+1)S_n)^{1/(n-1)}/H_{\min}^{(n-2)/(n-1)}(H_{\max}/H_{\min})^{-(n-2)/(n-1)} \\
&< C_{\min}^{\ell+1,\infty}.
\end{aligned}$$

By taking $\underline{A} = C_{\max}^{\ell,\infty} + \gamma$ with a small $\gamma > 0$ so that $\overline{A} < C_{\min}^{\ell+1,\infty}$, it is easy to see that the functional $J_H$ does not have any critical point nor critical point at infinity between the levels $\underline{A}$ and $\overline{A}$. Hence the desired result follows from Lemma 4.1. □

**Proof of Theorem 1.4:** We Argue by contradiction. We assume that the functional $J_H$ does not have any critical point. In particular, $J_H$ has no critical points under the level $C_{\min}^{2,\infty}$. Under the assumption of Theorem 1.4, it follows from Proposition 4.2 (with $k = 1$) that $J^{C_{\max}^{1,\infty}+\gamma}$ is a contractible set, for a small constant $\gamma$. Furthermore it is a retract by deformation of $C_{\min}^{2,\infty}$.
Using now the characterization of critical points at infinity in Proposition 3.4 and Corollary 3.6, we obtain that those who are under the level $C_{\min}^{2,\infty}$ are in one to one correspondence with critical points of $H$ in $\mathcal{H}$. Then we derive from Lemma 3.7 and the Euler-Poincaré theorem that:

$$1 = \chi(J^{C_{\min}^{2,\infty}+\gamma}) = \sum_{y \in \mathcal{H}} (-1)^{n-1-ind(H,y)}$$

which contradicts the assumption 2 of Theorem 1.4. Hence $J_H$ has at least one critical point. □

**Proof of Theorem 1.2:** We first point out that, under the assumption of the theorem, if $A_1 \neq 1$,

the existence of at least one solution to the problem (1.1) follows from Theorem 1.4. Hence we will assume that $A_1 = 1$ and we notice that the number $N := \sharp \mathcal{H}$ has to be odd. In fact, if $A_1 = 1$ we get $p - q = 1$ where $p$ and $q$ are respectively the number of critical points in $\mathcal{H}$ with even and odd Morse index. On the other hand, we have $p + q = N$. Then we obtain $2p = N + 1$ and the claim follows. We will come back later to this observation.

Assuming that $J_H$ does not have any critical point under the level $C_{\min}^{3,\infty}$, we derive, using Proposition 4.2 (with $k = 2$), that the level sets $J^{C_{\max}^{1,\infty}+\gamma}$ and $J^{C_{\max}^{2,\infty}+\gamma}$ are contractible sets. Then it follows from the properties of the Euler-Characteristic, see Proposition 5.7, pp.105 in [18], that

$$1 = \chi(J^{C_{\max}^{2,\infty}+\gamma}) = \chi(J^{C_{\max}^{2,\infty}+\gamma}, J^{C_{\max}^{1,\infty}+\gamma}) + \chi(J^{C_{\max}^{1,\infty}+\gamma}).$$

That is $\chi(J^{C_{\max}^{2,\infty}+\gamma}, J^{C_{\max}^{1,\infty}+\gamma}) = 0$. Moreover it follows from Proposition 3.4 and Corollary 3.6 that the critical points at infinity between these two levels are $(y_i, y_j)_\infty$ with $y_i \neq y_j \in \mathcal{H}$. Thus, it follows from Lemma 3.7 and the Euler-Poincaré theorem that

$$\sum_{y_i \neq y_j \in \mathcal{H}} (-1)^{1+\imath(y_i)+\imath(y_j)} = 0 \qquad (4.1)$$

where $\imath(y_k) := n - 1 - ind(H, y_k)$.

Recall that we are in the case of $A_1 = 1$. Then the number $N$ has to be odd, say $N := 2k + 1$ (with $k \in \mathbb{N}_0$) and there are $k$ odd numbers $\imath(y_j)$'s and $k+1$ even numbers $\imath(y_j)$'s. We claim that, for each $k \geq 0$, it holds

$$A_2 := \sum_{i<j} (-1)^{\imath(y_i)+\imath(y_j)} = -k. \qquad (4.2)$$

Indeed, observe that the value of $A_2$ is the sum of $+1$ and $-1$. To get $-1$, $\imath(y_i)$ and $\imath(y_j)$ have to be of different parity. However, to get $+1$, $\imath(y_i)$ and $\imath(y_j)$ have to be of the same parity. Hence

- For $k = 0$, we have only one point $y$ with an even $\imath(y)$. Thus $A_2 = 0$.

- For $k = 1$, we have two points $y_0$ and $y_2$ with even $\imath(y_k)$ and one point $y_1$ with an odd $\imath(y_1)$. Thus, $A_2 = 1 - 2 = -1$.

- For $k \geq 2$, there exist $k+1$ even numbers $\imath(y_j)$ and $k$ odd numbers $\imath(y_j)$. Thus, it holds

$$A_2 = \binom{k+1}{2} + \binom{k}{2} - \binom{k+1}{1}\binom{k}{1} = \frac{1}{2}(k+1)k + \frac{1}{2}k(k-1) - (k+1)k = -k.$$

Now, (4.1) and (4.2) imply that $k = 0$ and we get $\sharp \mathcal{H} = 1$. This leads to a contradiction with the assumption that $\sharp \mathcal{H} \geq 2$. Thereby the proof of the theorem is completed. □

**Acknowledgment:**
The author would like to thank Professor Mohameden Ahmedou and Mohamed Ben Ayed for having brought to his attention the problem addressed in the paper and giving him a lot of encouragement and suggestions throughout the work.